\documentclass[12pt,centertags]{amsart}
\usepackage{amsmath,amstext,amsthm,a4,amssymb,amscd}
\usepackage[mathscr]{eucal}
\usepackage{mathrsfs}
\usepackage{epsf}

\usepackage{amsmath,amstext,amsthm,a4,amssymb,amscd,mathrsfs}
\usepackage{mathtools}
\usepackage[mathscr]{eucal}
\usepackage{mathrsfs}
\usepackage{epsf}
\numberwithin{equation}{section}

\newcommand{\field}[1]{\mathbb{#1}}
\newcommand{\Z}{\field{Z}}

\newcommand{\spec}{\mathrm{spec}}
\newcommand{\pa}{\partial}
\newcommand{\bpa}{\bar{\partial}}
\newcommand{\bline}[1]{\overline{#1}}

\DeclareMathOperator{\Hom}{Hom}

\newtheorem{thm}{Theorem}[section]
\newtheorem{lemma}[thm]{Lemma}
\newtheorem{prop}[thm]{Proposition}
\newtheorem{cor}[thm]{Corollary}
\theoremstyle{definition}
\newtheorem{rem}[thm]{Remark}
\theoremstyle{definition}

\theoremstyle{definition}
\newtheorem{defn}[thm]{Definition}

\newcommand{\be}{\begin{eqnarray}}
\newcommand{\ee}{\end{eqnarray}}

\newcommand{\comment}[1]{}
\begin{document}

\title{A note on equivariant twisted holomorphic analytic torsion}
 
\author{Minxiang Yuan}

\address{Chern Institute of Mathematics \& LPMC, Nankai
University, Tianjin 300071, P.R. China}
\email{2120240021@mail.nankai.edu.cn}

\begin{abstract}
In this note, we introduce an equivariant version of twisted holomorphic analytic torsion for the twisted Dolbeault complex on a complex manifold equipped with an action of a compact Lie group. Our construction is motivated by Ray and Singer's holomorphic analytic torsion, by the equivariant analytic torsion developed for elliptic complexes with group actions, notably in the work of Donnelly and later in the holomorphic setting by Bismut, and by Mathai and Wu's twisted analytic torsion for the twisted de Rham complex. Using pseudodifferential operators and the equivariant noncommutative residue, we define the equivariant twisted holomorphic analytic torsion. We also state its variation formulas under changes of the metric and of the twisting element, and present its functorial properties.
\end{abstract}
\maketitle

\section{Introduction}
Let $X$ be a compact complex manifold, let $E$ be a holomorphic vector bundle over $X$, and let
$H \in \Omega^{\bar{1}}(X)$ satisfy $\bar{\partial}H=0$.
Here $\Omega^{\bar{0}}(X)$, respectively $\Omega^{\bar{0}}(X,E)$, denotes the direct sum of all, respectively all $E$-valued, $(0,q)$-forms on $X$ with $q$ even; similarly,
$\Omega^{\bar{1}}(X)$, respectively $\Omega^{\bar{1}}(X,E)$, denotes the direct sum of all, respectively all $E$-valued, $(0,q)$-forms on $X$ with $q$ odd. In this situation, one can define the twisted Dolbeault cohomology \cite[Section 7]{MathaiWuZ2gradedtorsion}:
\begin{equation}
    H^{\bline{k}}(X,E,H): =
    \frac{\mathrm{Ker}(\bpa_{\bline{k}}^{E, H}:\Omega^{\bline{k}}(X, E)\to \Omega^{\bline{k+1}}(X, E))}
    {\mathrm{Im}(\bpa_{\bline{k}-1}^{E, H}:\Omega^{\bline{k-1}}(X, E)\to \Omega^{\bline{k}}(X, E))},
    \quad k=0,1,
\end{equation}
where $\bpa_{\bline{k}}^{E,H}:=\bpa_{\bline{k}}^E+H\wedge\cdot$ is the twisted Dolbeault operator, and we use $\bline{k}$ to denote the remainder of $k$ modulo $2$. In general, the twisted differential preserves only the $\Z/2\Z$ grading. We call $H$ the \textbf{twisting element} of the twisted Dolbeault complex. 

Moreover, since $X$ is compact, one can prove that the twisted Dolbeault cohomology groups $H^{\bline{k}}(X,E,H)$, $k=0,1$, are finite-dimensional.

Now assume that a compact Lie group $G$ acts holomorphically on $X$ and $E$, and $H$ is $G$-invariant. Consider a Riemannian metric $g^{TX}$ on $X$ which is compatible with
the complex structure $J$ and is $G$-invariant. That is, for any vectors
$u,v \in TX$, we have
\begin{equation*}
    g^{TX}(Ju,Jv)=g^{TX}(u,v),
\end{equation*}
and the action of every element of $G$ on $X$ is by isometries. Moreover, let
$h^E$ be a $G$-invariant Hermitian metric on $E$. Then, on
$\Omega^{0,q}(X,E)$(so on $\Omega^{\bline{k}}(X,E)$), we have a $G$-invariant Hermitian inner product $\langle \cdot,\cdot \rangle$ such that
for $u=\alpha f$ and $v=\beta s$, where
$\alpha,\beta \in \Omega^{0,q}(X)$ and $f,s \in \Gamma(X,E)$,
\begin{equation}
    \langle u,v \rangle
    =
    \int_X (\alpha \wedge *\beta) h^E(f,s),
\end{equation}
where the operator $*$ is the Hodge star operator defined by $g^{TX}$.
The $G$-invariance means that
\begin{equation*}
    \langle g_*u,g_*v\rangle=\langle u,v\rangle.
\end{equation*}
    
\begin{rem}\label{rem:hodge}
    Here the operator $*$ is obtained by taking the complex conjugate of the
    complex-linear extension of the real Hodge star operator. According to
    \cite[Definition 3.1.3]{MR2093043}, the complex-linear extension of the Hodge star operator
    maps $(p,q)$-forms to $(n-q,n-p)$-forms. Therefore, the operator $*$ here
    maps $\Omega^{p,q}(X,E)$ to $\Omega^{n-p,n-q}(X,E)$, and is a
    conjugate-linear operator.
\end{rem}

By the Hodge theorem for the twisted elliptic complex, the twisted Dolbeault cohomology group $H^{\bline{k}}(X,E,H)$ can be viewed as a finite-dimensional subspace  of $\Omega^{\bline{k}}(X,E)$. Equip the cohomology with the inner
product induced by its harmonic
representatives. Following the construction in \cite{MR1262703}, we can construct an equivariant metric $\|\cdot\|$ on
\begin{equation}
      \det H^*(X,E,H)
    :=
    \det H^{\bline{0}}(X,E,H)
    \otimes
    \det H^{\bline{1}}(X,E,H)^{-1}.
\end{equation}
For details, see Definition \ref{def:equivariantmetric} in Section 3 or \cite{Bismut2}.

On the other hand, consider the adjoint operator $(\bpa_{\bline{k}}^{E,H})^*$ of $\bpa_{\bline{k}}^{E,H}$, as well as the corresponding twisted Kodaira--Laplace operator $\square_{\bline{k}}^{E,H}$,
\begin{equation}
    \square_{\bline{k}}^{E,H} = (\bpa_{\bline{k}}^{E,H})^* \bpa_{\bline{k}}^{E,H}+\bpa_{\bline{k-1}}^{E,H}(\bpa_{\bline{k-1}}^{E,H})^*, k=0,1.
\end{equation}

The operator $\square_{\bline{k}}^{E,H}$ is also a $G$-invariant elliptic operator and yields the Hodge decomposition
\begin{equation}
    \Omega^{\bline{k}}(X,E)
    \cong
    \ker \square_{\bline{k}}^{E,H}
    \oplus
    \mathrm{Im}(\bpa_{\bline{k-1}}^{E,H})
    \oplus
    \mathrm{Im}((\bpa_{\bline{k}}^{E,H})^*).
\end{equation}

Now we can define the equivariant $\zeta$-function of $(\bpa_{\bline{k}}^{E,H})^*\bpa_{\bline{k}}^{E,H}$ by
\begin{equation}
    \zeta_g\bigl(s,(\bpa_{\bline{k}}^{E,H})^*\bpa_{\bline{k}}^{E,H}\bigr)
    :=
    \mathrm{Tr}'\bigl(gP_{\overline{k}}(\square_{\bline{k}}^{E,H})^{-s}\bigr),
\end{equation}
where $g\in G$, $P_{\overline{k}}$ is the projection from $\Omega^{\bline{k}}(X,E)$ onto $\mathrm{Im}(\bpa_{\bline{k}}^{E,H})^*$, $(\square_k^{E,H})^{-s}$ is the complex power of the Laplacian \cite[Chapter 2]{Pseudodifferential}, and $\mathrm{Tr}'$ denotes the trace with the kernel of $\square_k^{E,H}$ removed. If $\mathrm{Re}(s)$ is sufficiently large, this function is well-defined and holomorphic in $s$.

By \cite[Proposition 2.1]{equnonresidue}, when $G$ is finite, the function
$\zeta_g\bigl(s,(\bpa_{\bline{k}}^{E,H})^*\bpa_{\bline{k}}^{E,H}\bigr)$
admits a meromorphic continuation to the entire complex plane, with at most a simple pole at $s=0$. Acutually the same conclusion holds when $G$ is a compact
Lie group, since the stationary-phase argument used in the proof remains
valid in this setting.

In this note, we prove the following result:
\begin{thm}\label{thm:main}
    The meromorphic continuation of the equivariant $\zeta$-function
    \begin{equation*}
        \zeta_g\bigl(s,(\bpa_{\bline{k}}^{E,H})^*
        \bpa_{\bline{k}}^{E,H}\bigr)
    \end{equation*}
    to the complex plane, still denoted by the same symbol, has residue zero at $s=0$. Equivalently, it is holomorphic at $s=0$.
\end{thm}

Consequently, we can define the equivariant $\zeta$-regularized determinant of the operator $(\bpa_{\bline{k}}^{E,H})^*\bpa_{\bline{k}}^{E,H}$ by
\begin{equation}
    \mathrm{Det}'_g\bigl[
    (\bpa_{\bline{k}}^{E,H})^*
    \bpa_{\bline{k}}^{E,H}
    \bigr]
    :=
    \exp\left(
    -\left.
    \frac{d}{ds}
    \zeta_g\bigl(s,(\bpa_{\bline{k}}^{E,H})^*
    \bpa_{\bline{k}}^{E,H}\bigr)
    \right|_{s=0}
    \right).
\end{equation}

By computing the derivative of $\zeta_g\bigl(s,(\bpa_{\bline{k}}^{E,H})^*
    \bpa_{\bline{k}}^{E,H}\bigr)$ at $s=0$, we can define the equivariant twisted holomorphic analytic torsion $\tau_g(X,E,H)$. It is given by the following expression, see Definition \ref{def:equtorsion}:
\begin{equation}
    \tau_g(X,E,H)
    :=
    \mathrm{Det}'_g(\bpa_{\bline{1}}^*\bpa_{\bline{1}})^{1/2}
    \mathrm{Det}'_g(\bpa_{\bline{0}}^*\bpa_{\bline{0}})^{-1/2}.
\end{equation}

We call this expression the \textbf{twisted holomorphic analytic torsion under the action of a compact group}, or \textbf{equivariant twisted holomorphic analytic torsion}.

We can also define \textbf{the equivariant twisted metric $\|\cdot\|_{E,H}(g)$} on $\det H^*(X,E,H)$ by
\begin{equation}
    \|\cdot\|_{E,H}(g):= \tau_g(X,E,H)\|\cdot\|,
\end{equation}

There are also some variation formulas under variations of $g^{TX}$, $h^E$, and $H$, and some important functorial properties for $\|\cdot\|_{E,H}(g)$, including the multiplicative decomposition for direct sums of vector bundles, the Serre duality relation, and the K\"unneth formula for products of manifolds.

The rest of this note is organized as follows.  In Section 2, we will give a proof of Theorem \ref{thm:main}. In Section 3, we will give the explicit definition of  the equivariant twisted holomorphic analytic torsion $\tau_g(X,E,H)$ and equivariant twisted metric $\|\cdot\|_{E,H}(g)$. In Section 4, we will present the explicit formulas for  variation formulas and functorial properties of $\|\cdot\|_{E,H}(g)$, which can be obtained by by adapting the arguments of  \cite[Section 3, Section 4]{MathaiWutwistedtorsion}.

\section{Proof of Theorem \ref{thm:main}}
\subsection{Equivariant noncommutative residue}
We briefly introduce the equivariant noncommutative residue theory, which plays a key role in the proof. For more details, we refer to \cite[Section 2]{equnonresidue}.

 Let $M$ be a closed manifold. Fix a first-order elliptic pseudodifferential operator $D \in \Psi^1(M)$ and assume that $D$ is positive. Let a compact group $G$ act smoothly on $M$, and assume that $D$ is a $G$-invariant operator. That is, $gD=Dg$ holds for every $g \in G$. Consider the function $z \mapsto \mathrm{Tr}(D^{-z}Ag)$. When $\mathrm{Re}(z)>n+\mathrm{ord}A$, $D^{-z}A$ is a trace-class operator (trace-class operators form an ideal in the space of bounded operators), and hence the function $z \mapsto \mathrm{Tr}(D^{-z}Ag)$ is holomorphic on the half-plane $\mathrm{Re}(z)>n+\mathrm{ord}A$.

Denote the function $z \mapsto \mathrm{Tr}(D^{-z}Ag)$ by $\zeta_{g,A}(z)$. We have the following theorem \cite[Theorem 2.1]{equnonresidue}.

\begin{thm}\label{thm:equivaresidue}
    Let $k_g=\dim M^g$. Let $D$ be a positive first-order elliptic differential operator which is $G$-invariant, and let $A$ be a classical pseudodifferential operator. Then the function $\zeta_{g,A}(z)$ is holomorphic on $\mathrm{Re}(z)>d$, where $d=k_g+\mathrm{ord}(A)$, and admits a meromorphic continuation to the whole complex plane. Moreover, the only possible poles of $\zeta_{g,A}(z)$ are simple poles at $z=d,d-1,d-2,\dots$.
\end{thm}
\begin{rem}\label{rem:holofamily} \cite[Remark 2.2]{equnonresidue}
   if $D$ has order $m$, then the theorem still holds, but the possible simple poles are $z=\frac{d}{m},\frac{d-1}{m},\frac{d-2}{m},\dots$.
\end{rem}

\begin{rem}\cite[the discussion following Proposition 2.1]{equnonresidue}
    In the proof of proposition \ref{thm:equivaresidue}, the argument uses only the fact that $D^{-z}$ is a holomorphic family and has order $-z$. Therefore, $D^{-z}$ can be replaced by any holomorphic family $R(z)$ whose order is $-z$. 
\end{rem}

By the preceding theorem, we can define the equivariant noncommutative residue of the operator $A$.

Although the function $\zeta_{g,A}=\mathrm{Tr}(D^{-z}Ag)$ depends on the choice of the $G$-invariant operator $D$, its residue at $z=0$ does not.
\begin{prop}\cite[p.~713, top paragraph]{equnonresidue}
    Let $D$ and $D'$ both be positive first-order $G$-equivariant elliptic operators. Let $\zeta_{g,A}=\mathrm{Tr}(gD^{-z}A)$ and $\zeta'_{g,A}=\mathrm{Tr}(gD'^{-z}A)$. Then $\mathrm{Res}_{z=0}\zeta_{g,A}=\mathrm{Res}_{z=0}\zeta'_{g,A}$.
\end{prop}

Thus we can define the equivariant noncommutative residue as follows.
\begin{defn}
    Let $D$ be a positive first-order $G$-equivariant elliptic operator. For any pseudodifferential operator $A$, we define the \textbf{equivariant noncommutative residue} $\mathrm{Res}_g(A)$ to be the residue at $z=0$ of the function $\zeta_{g,A}$.
\end{defn}

If $g=e$, then $\mathrm{Res}_g(A)$ reduces to the noncommutative residue. In addition, when $A$ is a smoothing operator, $\zeta_{g,A}$ is defined for every $s$, and hence $\mathrm{Res}_g(A)=0$. Therefore $\mathrm{Res}_g$ descends to a linear functional on $\Psi(M)/\Psi^{-\infty}(M)$.

By Theorem \ref{thm:equivaresidue}, we can also obtain the following corollary:
\begin{cor}\label{cor:residuevanish}
    If the action of $g$ on $M$ has only isolated fixed points, and if the action of $g$ on $S^*M$ has no fixed points, then $\mathrm{Res}_g(A)=0$.

    If $A$ is a smoothing operator, then $\mathrm{Res}_g(A)=0$.
\end{cor}
\begin{proof}
    According to the proof of Theorem \ref{thm:equivaresidue}, the expression for the residue is an integral over the cosphere bundle $S^*M^g$ of the fixed-point submanifold. But isolated fixed points mean that $S^*M^g$ is empty, and hence there is no residue.

    If $A$ is a smoothing operator, then $\zeta_{g,A}(z)$ is is entire in z, and hence $\mathrm{Res}_g(A)=0$.
\end{proof}

If we assume that $A$ is $G$-invariant, the residue also has the trace property: it vanishes on commutators.
\begin{prop}\label{prop:residuecommute}
    Let $A$ and $B$ be $G$-invariant pseudodifferential operators. Then $\mathrm{Res}_g(AB)=\mathrm{Res}_g(BA)$.
\end{prop}
\begin{proof}
    The method used here is the same as \cite[Corollary 2.3]{equnonresidue}. Let $D$ be a positive first-order $G$-invariant elliptic operator. Then
    \begin{equation}
    \begin{aligned}
        \mathrm{Res}_g([A,B])&=\mathrm{Res}_{z=0}\mathrm{Tr}(gD^{-z}[A,B])\\
        &=\mathrm{Res}_{z=0}[\mathrm{Tr}(gD^{-z}AB)-\mathrm{Tr}(gD^{-z}BA)]\\
        &=\mathrm{Res}_{z=0}[\mathrm{Tr}(gBD^{-z}A)-\mathrm{Tr}(gBA D^{-z})]\\
        &=\mathrm{Res}_{z=0}[\mathrm{Tr}(gB[D^{-z},A])].
    \end{aligned}
    \end{equation}
    
    At $z=0$, we have $[D^{-z},A]=0$, and hence $[D^{-z},A]=zR(z)$ for some holomorphic family of operators $R(z)$. Since $\mathrm{Tr}(gBR(z))$ has at most a simple pole at $z=0$, we have $\mathrm{Res}_{z=0}[\mathrm{Tr}(gB[D^{-z},A])]=0$. Therefore $\mathrm{Res}_g(AB)=\mathrm{Res}_g(BA)$.
\end{proof}
\begin{rem}
    Proposition \ref{prop:residuecommute} is in fact a special case of  \cite[Corollary 2.3]{equnonresidue}.
\end{rem}
\begin{rem}
    In the definition of the equivariant noncommutative residue, we always assume that the order of $D$ is $1$. If the order of $D$ is $m$, then by considering $D^{1/m}$ it is not hard to obtain
    \begin{equation}
        \mathrm{Res}_{z=0}\mathrm{Tr}(gD^{-z}A)=\frac{1}{m}\mathrm{Res}_g(A).
    \end{equation}

    If the operator $D$ is only nonnegative, then one may consider  the trace on the orthogonal complement of the kernel of $D$: 
    \begin{equation}
        \mathrm{Res}_{z=0}\mathrm{Tr}'(gD^{-z}A)=\mathrm{Res}_g(A).
    \end{equation}
    Since the Kodaira--Laplace operator considered later is only nonnegative, we adopt this method.
\end{rem}

Now consider the expression $ \mathrm{Tr}'(gP_{\overline{k}}(\square_{\bline{k}})^{-s})$. In order to obtain the regularity of $\mathrm{Tr}'(gP_{\overline{k}}(\square_{\bline{k}})^{-s})$ at $s=0$, we only need to prove $\mathrm{Res}_g(P_{\bline{k}})=0$. 

\subsection{When $g\neq e$ has only isolated fixed points}
Suppose that $g \neq e \in G$ has only isolated fixed points on $X$. In this situation, by \cite[Remark 2.5]{equnonresidue}, we know that the equivariant noncommutative residue vanishes. That is,
\begin{equation}
    \mathrm{Res}_g (P_{\bline{k}})=0, \quad k=0,1
\end{equation}
\subsection{General case}
Now we put no additional 
hypothesis on the action of $G$. Suppose $H$ can be written as
\begin{equation}
    H= \sum_{i  \geq 0}H_{(2i+1)}, \text{ where } H_{(2i+1)}\in \Omega^{0,2i+1}(X).
\end{equation}

We first state two lemmas. Since the proof of the first lemma follows the same argument as that of \cite[Proposition 1.1]{MathaiWutwistedtorsion}, we omit it.
\begin{lemma}\label{lem:global}
    Let $\lambda\geq 0$ be a non-negative real number. Define operator  $c_\lambda$ by multiplication by $\lambda^{[i/2]}$ on $\Omega^{0,i}(X,E)$.(let $\lambda^0 =1$ when $\lambda=0$) Then the following
identities hold:
    \begin{equation*}
        c_\lambda \circ \bpa^{E,H}_{\bline{k}} =\lambda^{k}\bpa^{E,{c_\lambda(H)}}_{\bline{k}} \circ c_\lambda,k=0,1.
    \end{equation*}

    Therefore, $c_\lambda$ induces an isomorphism between $H^{\bline{k}}(X,E,H)$ and $H^{\bline{k}}(X,E,c_\lambda(H))$ when $\lambda>0$.
\end{lemma}
\begin{lemma}
    Let $P_t$ be a smooth family of $G$-invariant 0-order pseudodifferential operators, with $P_t^2=P_t$ for all $t$. Fix $g\in G$. Then the equivariant residue $\mathrm{Res}_g(P_t)$ is independent of $t$.
\end{lemma}
\begin{proof}
    Let $\dot{P}_t$ be the derivative of $P_t$ with respect to $t$. Then we have $\dot{P}_tP_t+P_t\dot{P}_t=\dot{P}_t$. Therefore,
    \begin{equation*}
         \mathrm{Res}_g (  \dot{P_t})  = \mathrm{Res}_g(  P_t\dot{P_t}) + \mathrm{Res}_g(  \dot{P_t}P_t). 
    \end{equation*}
    As $P_t$ is $G$-invariant, $\dot{P_t}$ is $G$-invariant. By Proposition \ref{prop:residuecommute}, we have: 
    \begin{equation}
        \begin{aligned}
            \mathrm{Res}_g(  P_t\dot{P_t}) =\mathrm{Res}_g(  P_t^2\dot{P_t}) =\mathrm{Res}_g(  P_t\dot{P_t}P_t) =0. \\
        \end{aligned}
    \end{equation}

    Similarly, we have $\mathrm{Res}_g(  \dot{P_t}P_t) =0$. Therefore, $\mathrm{Res}_g \dot{P_t}=0$. 
\end{proof}

Now we can give a proof of Theorem \ref{thm:main}.

Deform $H$ by setting $H_\lambda:=c_\lambda(H),\lambda\in [0,1]$. So $H_0=H_{(1)}$, $H_1=H$. Then we have a family of twisted Dolbeault operators $\bpa_{\bline{k}}^{E,H_\lambda}$, and the corresponding projection $P_{\overline{k},\lambda}$ from $\Omega^{\bline{k}}(X,E)$ onto $\mathrm{Im}(\bpa_{\bline{k}}^{E,H_\lambda})^*$. For simplicity, we set $\bpa_{\bline{k},\lambda}=\bpa_{\bline{k}}^{E,H_\lambda}$, $(\bpa_{\bline{k},\lambda})^*=(\bpa_{\bline{k}}^{E,H_\lambda})^*$, and $\square_{\bline{k},\lambda}=\square_{\bline{k}}^{E,H_\lambda}$. Notice that $H_\lambda$, $\bpa_{\bline{k},\lambda}$, $(\bpa_{\bline{k},\lambda})^*$, $\square_{\bline{k},\lambda}$, and $P_{\overline{k},\lambda}$ are all $G$-invariant.

If $\mu$ is an eigenvalue of $\square_{\bline{k},\lambda}$, then let $\Omega^{\bline{k}}_{\{\mu\}}(X,E)_\lambda$ be the corresponding eigenspace, which is a finite-dimensional representation of $G$ by ellipticity of $\square_{\bline{k},\lambda}$. Moreover, projection from $\Omega^{\bline{k}}(X,E)$ onto $\Omega^{\bline{k}}_{\{\mu\}}(X,E)_\lambda$ is a $G$-invariant smoothing operator. 

Suppose $a>0$ is a positive number. Set
\begin{equation}
    \Omega^{\bar{k}}_{[0,a)}(X,E)_\lambda=\bigoplus_{\mu \in [0,a)} \Omega^{\bline{k}}_{\{\mu\}}(X,E)_\lambda.
\end{equation} 

Therefore, $\Omega^{\bline{k}}_{[0,a)}(X,E)_\lambda$ is a finite-dimensional representation of $G$, and the projection from $\Omega^{\bline{k}}(X,E)$ onto $\Omega^{\bline{k}}_{[0,a)}(X,E)_\lambda$ is a $G$-invariant smoothing operator, which we denote by $Q_{\overline{k},\lambda}^a$. Let $\Omega^{\bline{k}}_{[a,+\infty)}(X,E)_\lambda$ be the orthogonal complement of $\Omega^{\bline{k}}_{[0,a)}(X,E)_\lambda$ in $\Omega^{\bline{k}}(X,E)$. Then we have a $G$-invariant orthogonal decomposition
\begin{equation}
    \Omega^{\bline{k}}(X,E)
    =
    \Omega^{\bline{k}}_{[0,a)}(X,E)_\lambda
    \oplus
    \Omega^{\bline{k}}_{[a,+\infty)}(X,E)_\lambda.
\end{equation}

Denote the projection from $\Omega^{\bline{k}}(X,E)$ onto $\Omega^{\bline{k}}_{[a,\infty)}(X,E)_\lambda$ by $(Q_{\overline{k},\lambda}^a)^{\perp}$. Then we have 
\begin{equation*}
   P_{\overline{k},\lambda}=P_{\overline{k},\lambda}Q_{\overline{k},\lambda}^a+P_{\overline{k},\lambda}(Q_{\overline{k},\lambda}^a)^{\perp}, 
\end{equation*}
and $P_{\bline{k},\lambda}(\Omega^{\bline{k}}_{[0,a)}(X,E)_\lambda)\subset \Omega^{\bline{k}}_{[0,a)}(X,E)_\lambda$,  $P_{\bline{k},\lambda}(\Omega^{\bline{k}}_{[a,\infty)}(X,E)_\lambda)\subset \Omega^{\bline{k}}_{[a,\infty)}(X,E)_\lambda$. This gives a decomposition of $P_{\bline{k},\lambda}$.

It is clear that $P_{\bline{k},\lambda}Q_{\bline{k},\lambda}^a$ is also a smoothing operator. By Corollary \ref{cor:residuevanish}, its equivariant residue is $0$. Therefore,
\begin{equation*}
    \mathrm{Res}_g P_{\bline{k},\lambda}= \mathrm{Res}_g P_{\bline{k},\lambda}(Q_{\bline{k},\lambda}^a)^{\perp}, \forall \lambda \in [0,1].
\end{equation*}
For simplicity, now we write $P_{\bline{k},\lambda}(Q_{\bline{k},\lambda}^a)^{\perp}$ as $P_{\bline{k},\lambda}^a$.

We then state three propositions which are necessary for the proof.
\begin{prop}\label{prop:spectrum}
    Let $\spec(\square_{\bline{k},\lambda})=\{\mu_1(\lambda),\mu_2(\lambda),\dots\}$,where $\mu_i(\lambda)\leq \mu_{i+1}(\lambda),\forall i$. Then for any $i$, $\mu_i(\lambda)$ is a continuous function of $\lambda$.
\end{prop}
\begin{proof}
Consider the min-max principle for self-adjoint operators \cite[Theorem 4.10]{Teschl}. We have
$$
\mu_i(\lambda)=\sup_{v_1,\dots,v_{i-1}} \inf_{v \in U(v_1,\dots,v_{i-1})} 
\langle \square_{\bline{k},\lambda} v, v\rangle,
$$
where
$$
U(v_1,\dots,v_{i-1})
=
\{v \in H^{2,\bline{k}}(X,E) : v \perp v_1,\dots,v_{i-1},\|v\|_{L^2}=1\},
$$
and $H^{2,\bline{k}}(X,E)$ denotes the Sobolev space of order $2$ associated with $\Omega^{\bline{k}}(X,E)$.

Let
$$
D_{\overline{k},\lambda}
=
\bpa_{\bline{k},\lambda}
+
\bpa_{\bline{k+1},\lambda}^*.
$$
It is easy to verify $D_{\bline{k+1},\lambda}\circ D_{\bline{k},\lambda}=\square_{\bline{k},\lambda}$. Moreover, by
$$
\langle D_{\overline{k},\lambda} \cdot, \cdot \rangle
=
\langle \cdot, D_{\overline{k+1},\lambda} \cdot \rangle,
$$
we have
\begin{equation*}
    \langle \square_{\bline{k},\lambda}v,v\rangle
    =
    \langle D_{\bline{k+1},\lambda}\circ D_{\overline{k},\lambda}v,v \rangle
    =
    \langle D_{\bline{k},\lambda}v,D_{\bline{k},\lambda}v \rangle
    =
    \|D_{\bline{k},\lambda} v\|_{L^2}^2.
\end{equation*}

Let $\lambda_0 \in [0,1]$ be fixed. We claim that there exists a constant $C_1>0$ such that
\begin{equation*}
    \bigg | \langle (\square_{\bline{k},\lambda}-\square_{\bline{k},\lambda_0})v, v\rangle\bigg |
    \leq
    C_1 |\lambda-\lambda_0| \cdot \|v\|_{H^1}^2
\end{equation*}
holds for every $\lambda \in [0,1]$.

We now prove this claim. Consider
\begin{equation}\label{eq:minus}
    \begin{aligned}
       \bigg|\langle (\square_{\bline{k},\lambda} - \square_{\bline{k},\lambda_0}) v, v \rangle\bigg|
       &=
       \bigg|
       \langle D_{\bline{k},\lambda}v,D_{\bline{k},\lambda}v\rangle
       -
       \langle D_{\bline{k},\lambda_0}v,D_{\bline{k},\lambda_0}v \rangle
       \bigg| \\
       &=
       \bigg |
       \langle D_{\bline{k},\lambda}v,[D_{\bline{k},\lambda}-D_{\bline{k},\lambda_0}]v\rangle
       +
       \langle [D_{\bline{k},\lambda}-D_{\bline{k},\lambda_0}]v,D_{\bline{k},\lambda_0}v\rangle
       \bigg |.
    \end{aligned}
\end{equation}

We regard the Hermitian metric $h^E$ as a real-linear bundle map
$E \to E^*$ and use the same notation for its extension to
bundle-valued differential forms. Let
$\Gamma = h^E \circ * \colon \Omega^{p,q}(X,E)
\to \Omega^{n-p,n-q}(X,E^*)$, where the Hodge star operator $*$
is taken with the convention specified in Remark~\ref{rem:hodge}.
Following the argument in the proof of
\cite[Theorem 1.2]{MathaiWutwistedtorsion}, one verifies that
$\Gamma$ is an isometry and that
\[
    (\bpa_{\bline{k}}^{E,H})^*
    = (-1)^{k+1}\Gamma^{-1}
    \circ \bpa_{\bline{n+1-k}}^{E^*\otimes K_X,-H}
    \circ \Gamma.
\]

A direct computation gives
\begin{align*}
    D_{\bline{k},\lambda} - D_{\bline{k},\lambda_0}
    &=
    \bpa_{\bline{k},\lambda}
    - \bpa_{\bline{k},\lambda_0}
    + \bpa_{\bline{k+1},\lambda}^*
    - \bpa_{\bline{k+1},\lambda_0}^* \\
    &=
    (H_\lambda - H_{\lambda_0}) \wedge \cdot
    + (-1)^k \Gamma^{-1}
    \circ \bigl[(-H_\lambda + H_{\lambda_0}) \wedge \cdot\bigr]
    \circ \Gamma.
\end{align*}
Moreover, the identity
\[
    H_\lambda - H_{\lambda_0}
    =
    (\lambda - \lambda_0)
    \sum_{i \geq 1}
    \left(
        \sum_{j=0}^{i-1} \lambda^{i-1-j}\lambda_0^j
    \right)
    H_{2i+1}
\]
yields, for $|\lambda|, |\lambda_0| \leq 1$, the estimate
\[
    \|(H_\lambda - H_{\lambda_0})\wedge \cdot\|_{L^2}
    \leq
    |\lambda - \lambda_0|
    \sum_{i \geq 1} i \|H_{2i+1}\wedge \cdot\|_{L^2}.
\]

Since $\Gamma$ is an isometry and $\|H_{2i+1}\wedge \cdot\|_{L^2}$ is finite for each $i$, there exists $C_1'>0$ such that, for the fixed $\lambda_0$, one has
\begin{equation*}
    \|D_{\bline{k},\lambda}-D_{\bline{k},\lambda_0}\|_{H^1}
    \leq
    \|D_{\bline{k},\lambda}-D_{\bline{k},\lambda_0}\|_{L^2}
    \leq
    C_1' |\lambda-\lambda_0|,
\end{equation*}
where $\|D_{\bline{k},\lambda}-D_{\bline{k},\lambda_0}\|_{H^1}$ denotes the operator norm of this operator viewed as a map from $H^{1,\bline{k}}(X,E)$ to $L^2(X,E)$. Hence $\|D_{\bline{k},\lambda}\|_{H^1}$ is a continuous function of $\lambda$. Therefore, there exists $C_1''$ such that, for every $\lambda \in [0,1]$, we have
$$
\|D_{\bline{k},\lambda}\|_{H^1} \leq C_1''.
$$

Thus, by \eqref{eq:minus}, we obtain
\begin{equation*}
    \begin{aligned}
         \bigg|\langle (\square_{\bline{k},\lambda} - \square_{\bline{k},\lambda_0}) v, v \rangle\bigg|
         &\leq
         \left(
         \| D_{\bline{k},\lambda}v\|_{L^2}
         +
         \|D_{\bline{k},\lambda_0}v\|_{L^2}
         \right)
         \|(D_{\bline{k},\lambda}-D_{\bline{k},\lambda_0})v\|_{L^2} \\
         &\leq
         \left(
         \|D_{\bline{k},\lambda}\|_{H^1}
         +
         \|D_{\bline{k},\lambda_0}\|_{H^1}
         \right)
         \|(D_{\bline{k},\lambda}-D_{\bline{k},\lambda_0})\|_{H^1}
         \|v\|_{H^1}^2 \\
         &\leq
         C_1 |\lambda - \lambda_0| \|v\|_{H^1}^2.
    \end{aligned}
\end{equation*}

On the other hand, since $\square_{\bline{k},\lambda_0}$ is an elliptic operator, by the elliptic estimate \cite[Theorem 5.2]{MR1031992}, there exists a constant $c_0 > 0$, independent of $v$, such that
\begin{equation*}
    \|v\|_{H^1}^2
    \leq
    c_0
    \left(
    \|v\|_{L^2}^2
    +
    \|(\bpa_{\bline{k},\lambda_0}+\bpa_{\bline{k+1},\lambda_0}^*) v\|_{L^2}^2
    \right)
    =
    c_0
    \left(
    \|v\|_{L^2}^2
    +
    \langle \square_{\bline{k},\lambda_0} v, v \rangle
    \right).
\end{equation*}

Substituting this inequality into the preceding estimate, there exists a constant $C'>0$ such that, for all $\lambda,\lambda_0 \in [0,1]$, we have
\begin{equation*}
    |\langle \square_{\bline{k},\lambda} v, v \rangle
    -
    \langle \square_{\bline{k},\lambda_0} v, v \rangle|
    \le
    C'|\lambda-\lambda_0|
    \left(
    \langle \square_{\bline{k},\lambda_0} v, v \rangle
    +
    \|v\|_{L^2}^2
    \right).
\end{equation*}

Now, for $v \in U(v_1,\dots,v_{i-1})$, we have $\|v\|_{L^2} = 1$. Let
$\epsilon=|\lambda-\lambda_0|$, then the above inequality can be rewritten as
$$
(1-C'\epsilon)\langle \square_{\bline{k},\lambda_0} v, v \rangle
-
C'\epsilon
\le
\langle \square_{\bline{k},\lambda} v, v \rangle
\le
(1+C'\epsilon)\langle \square_{\bline{k},\lambda_0} v, v \rangle
+
C'\epsilon.
$$

Since taking infimum and supremum preserves the inequalities, we obtain
\begin{equation*}
    (1-C'\epsilon)\mu_i(\lambda_0) - C'\epsilon
    \le
    \mu_i(\lambda)
    \le
    (1+C'\epsilon)\mu_i(\lambda_0) + C'\epsilon.
\end{equation*}

Letting $\epsilon \to 0$, we get $\mu_i(\lambda)\to \mu_i(\lambda_0)$. Therefore, $\mu_i(\lambda)$ is continuous with respect to $\lambda$.
\end{proof}
\begin{prop}\label{prop:neighbourhood}
For any $a>0$ such that $a \notin \spec(\square_{\bline{k}},0)$, there exists a neighbourhood $U_a$ of $0$ such that $P_{\bline{k},\lambda}^a$ depends smoothly on $\lambda \in U_a$.
\end{prop}

\begin{proof}
    We only need to show that there exists a neighbourhood $U_a$ such that the dimension of $\Omega^{\bline{k}}_{[0,a)}(X,E,H)_\lambda$ remains unchanged when $\lambda$ varies in $U_a$. But this comes from the continuity of the spectrum of $\square_{\bline{k},\lambda}$(Proposition \ref{prop:spectrum}).
\end{proof}

\begin{prop}\label{prop:global}
    $P_{\bline{k},\lambda}$ depends smoothly on $\lambda$ when $\lambda>0$.
\end{prop}
\begin{proof}
    By Lemma \ref{lem:global}, the dimension of $H^{\bline{k}}(X,E,H)$ is independent of $\lambda$ when $\lambda$ varies in $(0,1]$. 
\end{proof}
With Proposition \ref{prop:neighbourhood} and Proposition \ref{prop:global}, if $a \notin \spec(\square_{\bline{k},0})$, then there exists $\lambda_0>0$ such that:
\begin{equation}
    \mathrm{Res}_g P_{\bline{k}}=\mathrm{Res}_g P_{\bline{k},1}=\mathrm{Res}_g P_{\bline{k},\lambda_0}=\mathrm{Res}_g P_{\bline{k},\lambda_0}^a=\mathrm{Res}_g P_{\bline{k},0}^a =\mathrm{Res}_g P_{\bline{k},0}
\end{equation}

When $\lambda=0$, we have $\bpa_{\bline{k},0}=\bpa_{\bline{k}}^{E,H_{(1)}}$, $\square_{\bline{k},0}=\square_{\bline{k}}^{E,H_{(1)}}$. Also, the twisted Dolbeault complex reduces and the differential has degree
one, so the complex is $\Z$-graded. We use $\square_{i}^{E,H_{(1)}}$ to denote the twisted Kodaira-Laplace operator restricted to $\Omega^{0,i}(X,E)$.

In this situation, we have:
\begin{equation*}
    \begin{aligned}
        \sum_{k=0,1}(-1)^k\zeta_g(s,\bpa_{\bline{k},0}^*\bpa_{\bline{k},0})&=\sum_{i=0}^n i(-1)^i \zeta_g(s,\square_{i}^{E,H_{(1)}})\\
        &=\sum_{i=0}^n i(-1)^i \mathrm{Tr}'(g(\square_{i}^{E,H_{(1)}})^{-s}).
     \end{aligned}
    \end{equation*}
    
But the regularity at $s=0$ of this expression is classical \cite[Definition 2.1]{Bismut2}. By a standard heat kernel analysis, we obtain that $\displaystyle\sum_{k=0,1}(-1)^k\zeta_g(s,\bpa_{\bline{k},0}^*\bpa_{\bline{k},0})$ is holomorphic at $s=0$, which implies that $\mathrm{Res}_g(P_{\bline{0},0}-P_{\bline{1},0})=\mathrm{Res}_g (P_{\bline{0}}-P_{\bline{1}})=0$.

    On the other hand, it is easy to verify that $\mathrm{Res}_g (P_{\bline{0}}+P_{\bline{1}})=\frac{1}{2}\mathrm{Res}_{s=0}\mathrm{Tr}'(g(\square_{\bline{0}})^{-s})=0$. Taking the two results above together, we obtain that
    \begin{equation*}
    \mathrm{Res}_g P_{\bline{0}}=\mathrm{Res}_g P_{\bline{1}}=0.
    \end{equation*}

 The proof of Theorem \ref{thm:main} now is completed.

\section{Definition of equivariant twisted holomorphic analytic torsion and the equivariant twisted metric}
\subsection{Equivariant metric}
Recall that $H^{\bline{k}}(X,E,H)$ is a finite-dimensional representation of $G$. According to the $G$-invariance of $\square^{E,H}$, we have a $G$-isomorphism $\ker \square_{\bline{k}}^{E,H} \cong H^{\bline{k}}(X,E,H), k=0,1$.  

\begin{defn}
    The \textbf{determinant space} of the twisted Dolbeault cohomology is defined as:
    \begin{equation*}
            \det H^{*}(X,E,H) = \bigotimes_{k=0,1} \det (H^{\bline{k}}(X,E,H))^{(-1)^k}.
    \end{equation*} 
\end{defn}

As $H^{\bline{k}}(X,E,H)$ is a finite-dimensional representation of $G$, $\det H^{*}(X,E,H)$ is also a finite-dimensional representation of $G$. 

Let $\hat{G}$ denote the set of irreducible representations of $G$. Then $H^{\bline{k}}(X,E,H)$ can be decomposed into its isotypical components:
\begin{equation}
    \begin{aligned}
        H^{\bline{k}}(X,E,H)
        &= \bigoplus_{(W,\rho_W)\in \hat{G}}
        \Hom_G(W,H^{\bline{k}}(X,E,H))\otimes W \\
        &= \bigoplus_{(W,\rho_W)\in \hat{G}}
        H^{\bline{k}}(X,E,H)_W\otimes W.
    \end{aligned}
\end{equation}
Here we use $H^{\bline{k}}(X,E,H)_W$ to denote the expression $ \Hom_G(W,H^{\bline{k}}(X,E,H)) \otimes W$

We can write $\det H^{\bline{k}}(X,E,H)$ as
\begin{equation*}
    \det H^{\bline{k}}(X,E,H)=\bigotimes_{(W,\rho_W)\in \hat{G}}\det H^{\bline{k}}(X,E,H)_W.
\end{equation*}

Consider the Hermitian metric on $H^{\bline{k}}(X,E,H)_W$ induced by the Hermitian metric on $H^{\bline{k}}(X,E,H)$. We can define the equivariant metirc on $\det H^{\bline{k}}(X,E,H)$.
\begin{defn}\cite[equation (0.6)]{Bismut2}
    Let $\| \cdot\|_{W,\bline{k}}$ denote the metric on $\det H^{\bline{k}}(X,E,H)_W$ induced by the Hermitian inner product on $H^{\bline{k}}(X,E,H)_W$. Consider the following expression:
    \begin{equation*}
         \log(||\cdot||_{\bline{k}})=\sum_{(W,\rho_W) \in \hat{G}} \log(||\cdot||_{W,\bline{k}}) \frac{\chi(W)}{\mathrm{rk}(W)},
    \end{equation*}
    where $\chi(W)$ denotes the character of $W$, $\mathrm{rk}(W)$ denotes the dimension of $W$.

    Then we define $||\cdot||_{\bline{k}}$  as the \textbf{equivariant metirc} of $\det H^{\bline{k}}(X,E,H)$.
\end{defn}

If we substitute $g \in G$ into $\|\cdot\|_{\bline{k}}$, we can get a metric on  $\det H^{\bline{k}}(X,E,H)$. We denote this metric as $\|\cdot\|_{\bline{k}}(g)$.

\begin{defn}\label{def:equivariantmetric}
    The \textbf{equivariant metric} $\|\cdot\| $ of $\det H^{*}(X,E,H)$ is defined as:
    \begin{equation*}
        \|\cdot\|:= \|\cdot\|_{\bline{0}}\otimes \|\cdot\|_{\bline{1}}^{-1}
    \end{equation*}
\end{defn}
If we substitute $g \in G$ into $\|\cdot\|$, we can get a metric on  $\det H^*(X,E,H)$. We denote this metric as $\|\cdot\|(g)$.

\subsection{Equivariant twisted holomorphic analytic torsion and equivariant twisted metric}

 By Theorem \ref{thm:main}, the equivariant $\zeta$-function $\zeta_g(s,\bpa_{\bline{k}}^* \bpa_{\bline{k}})=\mathrm{Tr}'(gP_{\bline{k}}(\square_{\bline{k}})^{-s})$ can be extended holomorphically at $s=0$. 
 
 Therefore, we can define the \textbf{$\zeta$-regularized determinant} of $\bpa_{\bline{k}}^* \bpa_{\bline{k}}$ as:
 \begin{equation}
    \mathrm{Det}'_g(\bpa_{\bline{k}}^* \bpa_{\bline{k}})= \exp(-\frac{d}{ds}\zeta_g(s,\bpa_{\bline{k}}^* \bpa_{\bline{k}})|_{s=0}).
 \end{equation} 

 Now we can define the equivariant twisted holomorphic analytic torsion and the equivariant twisted metric.
 \begin{defn}\label{def:equtorsion}
    Define the \textbf{equivariant twisted holomorphic analytic torsion} $\tau_g(X,E,H)$ by the following formula:
    \begin{equation}
        \tau_g(X,E,H):=\mathrm{Det}'_g(\bpa_{\bline{1}}^*\bpa_{\bline{1}})^{1/2}\mathrm{Det}'_g(\bpa_{\bline{0}}^*\bpa_{\bline{0}})^{-1/2}.
    \end{equation}
 \end{defn}
\begin{defn}
    Define the \textbf{equivariant twisted metric} $\|\cdot\|_{E,H}(g)$ by the following formula:
    \begin{equation}
        \|\cdot\|_{E,H}(g):= \tau_g(X,E,H) \cdot \|\cdot\|(g),
    \end{equation}
    where $\|\cdot\|(g)$ denotes the equivariant metric at $g \in G$. $\|\cdot\|_{E,H}(g)$ now is a metric on $\det H^*(X,E,H)$.
\end{defn}

\section{Variation formulas and functorial properties }
In this section, we state some properties of $\|\cdot\|_{E,H}(g)$ which can be obtained by the same method in \cite[Section 3, Section 4]{MathaiWutwistedtorsion}. In this section, we suppose $g \neq e$.

\subsection{Variation formula under changes of the
metrics}

Suppose $(g^{TX}_u,h^E_u)$ is a family of metrics on $(X,E)$ depending smoothly on $u$. Let $\alpha_u=*_u^{-1} \dfrac{\pa *_u}{\pa u}+(h_u^E)^{-1}\dfrac{\pa h_u^E}{\pa u} $, where $*_u$ is the Hodge star operator, $h^E$ is the real linear isomorphism from $E$ to $E^*$ induced by the Hermitian structure of $E$. Therefore, $\alpha_u$ is an operator from $\Omega^{\bline{k}}(X,E)$ to itself. 

\begin{prop}
    With the above notation, we have
    \begin{equation}
    \begin{aligned}
          &\frac{\pa}{\pa u}\log(\mathrm{Det}'_g(\bpa_{\bline{1}}^*\bpa_{\bline{1}})(g)\mathrm{Det}'_g(\bpa_{\bline{0}}^*\bpa_{\bline{0}})(g)^{-1})\\&=\sum_{k=0,1}(-1)^k (\mathrm{Lim}_{t \to 0^+}\mathrm{Tr}(g\alpha_u e^{-t\square_{\bline{k}}})-\mathrm{Tr}(g\alpha_u Q_{\bline{k},u})),
    \end{aligned}
    \end{equation}
    where $\mathrm{Lim}_{t \to 0^+}\mathrm{Tr}(g\alpha_u e^{-t\square_{\bline{k}}})$ denotes the constant term of the asymptotic expansion of $\mathrm{Tr}(g\alpha_u e^{-t\square_{\bline{k}}})$ at $t=0$, and $Q_{\bline{k},u}$ denotes the projection from $\Omega^{\bline{k}}(X,E)$ to $\ker \square_{\bline{k},u}$ by Hodge decomposition.
\end{prop}
 When $u$ varies, we can view $\|\cdot \|(g)$ as a family of metrics on $\det H^*(X,E,H)$ by making use of Hodge decomposition. More precisely, let $\|\cdot\|_{\bline{k},w}$ denote the equivariant metric of $\det H^{\bline{k}}(X,E,H)$ corresponding to parameter $w$. As $Q_{\bline{k},w}$ induces a $G$-isomorphism from $\ker \square_{\bline{k},u}$ to $\ker \square_{\bline{k},w}$, we can view $\|\cdot\|_{\bline{k},w}$ as a metric on $\det \ker \square_{\bline{k},u} \cong \det H^{\bline{k}}(X,E,H)$. Therefore, we can view $\|\cdot \|(g)$ as a family of metrics on $\det H^*(X,E,H)$. 
\begin{prop}
    Fix u as the base parameter. Then
    \begin{equation}
        \frac{\pa}{\pa w}\bigg|_u\log(\frac{\|\cdot\|^2_{w}(g)}{\|\cdot\|^2_{u}(g)})=\sum_{k=0,1}(-1)^k \mathrm{Tr}(g\alpha_u Q_{\bline{k},u}),
    \end{equation}
    where $Q_{\bline{k},u}$ denotes the projection from $\Omega^{\bline{k}}(X,E)$ to $\ker \square_{\bline{k},u}$ by Hodge decomposition.
\end{prop}

Now we can present the variation formula of metric.
\begin{thm}\label{thm:variationmetric}
    Let $X$ be a compact complex manifold, let $E$ be a holomorphic vector bundle
over $X$, and let $H$ be a $\bpa$-closed $(0,\bline{1})$-form. Suppose that
$G$ acts holomorphically on $X$. For a smooth one-parameter family of metrics
$(g_u^{TX},h_u^E)$, with $u$ as the parameter, we have
\begin{equation}
    \frac{\pa}{\pa u}\log[\|\cdot\|_{E,H}(g)^2]
    =
    \sum_{k=0,1}(-1)^k
    \mathrm{Lim}_{t \to 0^+}\mathrm{Tr}
    \bigl(g\alpha_u e^{-t\square_{\bline{k}}}\bigr),
\end{equation}
 where $\mathrm{Lim}_{t \to 0^+}\mathrm{Tr}(g\alpha_u e^{-t\square_{\bline{k}}})$ denotes the constant term of the asymptotic expansion of $\mathrm{Tr}(g\alpha_u e^{-t\square_{\bline{k}}})$ at $t=0$.
\end{thm}

\subsection{Variation formula for the twisting element $H$}
Suppose $H_v$ is a family of twisting elements on $(X,E)$ depending smoothly on $v$. We only consider the case when
the cohomology class $[H]$ remains fixed, which by \cite[Section 1.2]{MathaiWutwistedtorsion}, implies that the dimension of $H^{\bline{k}}(X,E,H)$ remains unchanged. We suppose that
\begin{equation}
    \frac{\pa H}{\pa v}=-\bpa B,
\end{equation}
where $B\in \Omega^{\bline{0}}(X)$. Set
\begin{equation*}
    \beta:=B\wedge \cdot:\Omega^{\bar{k}}(X,E)\to \Omega^{\bline{k}}(X,E),
\end{equation*}
then we have the following lemma.

\begin{lemma}
Under the above assumptions, for a smooth one-parameter family of $\bpa$-closed
$(0,\bline{1})$-forms $H_v$, we have
\begin{equation}
\begin{aligned}
&\frac{\pa}{\pa v}
\log\left(
\mathrm{Det}'_g(\bpa_{\bline{1}}^*\bpa_{\bline{1}})(g)
\mathrm{Det}'_g(\bpa_{\bline{0}}^*\bpa_{\bline{0}})(g)^{-1}
\right) \\
&=
\sum_{k=0,1}(-1)^k
\left[
\mathrm{Lim}_{t \to 0^+}\mathrm{Tr}(g\beta e^{-t\square_{\bline{k}}})
-\mathrm{Tr}(g\beta Q_{\bline{k},v})
\right] \\
&\quad+
\sum_{k=0,1}(-1)^k
\overline{
\left[
\mathrm{Lim}_{t \to 0^+}\mathrm{Tr}(g^*\beta e^{-t\square_{\bline{k}}})
-\mathrm{Tr}(g^*\beta Q_{\bline{k},v})
\right]
},
\end{aligned}
\end{equation}
where $g^*$ denotes the adjoint of the operator $g$, $\mathrm{Lim}_{t \to 0^+}\mathrm{Tr}(g\beta e^{-t\square_{\bline{k}}})$ denotes the constant term of the asymptotic expansion of $\mathrm{Tr}(g\beta e^{-t\square_{\bline{k}}})$ at $t=0$, and $Q_{\bline{k},v}$ denotes the projection from $\Omega^{\bline{k}}(X,E)$ to $\ker \square_{\bline{k},v}$ by Hodge decomposition.
\end{lemma}

We now consider the variation of the equivariant metric. By \cite[Section 1.2]{MathaiWutwistedtorsion}, as $H$ varies, we can use the isomorphism $\epsilon_B:=e^B \wedge \cdot$ to identify the corresponding cohomology groups. In this case, let $\bpa C = H_v-H_w$, then
\begin{equation*} 
    Q_{\overline{k},w}\circ \epsilon_{C}
\end{equation*}
is a $G$-equivariant isomorphism from
$\ker \square_{\overline{k},v}$ to $\ker \square_{\overline{k},w}$, whose inverse is given by
\begin{equation*}
Q_{\overline{k},v}\circ \epsilon_{-C}.
\end{equation*}
Taking $v$ as the base point, we can construct a family of metrics $\|\cdot\|_{\overline{k},w}$ on
\begin{equation*}
\det H^{\overline{k}}(X,E,H_v)\cong \det \ker \square_{\overline{k},v}
\end{equation*}
depending on $w$. Hence we obtain a family of metrics
$\|\cdot\|_w$ on $\det H^*(X,E,H_v)$, depending on $w$.

\begin{lemma}
Under the above assumptions, for a smooth one-parameter family of
$\overline{\partial}$-closed $(0,\overline{1})$-forms $H_v$, we have
\begin{equation}
    \left.\frac{\partial}{\partial w}\right|_{v}
    \log\left(
    \frac{\|\cdot\|^2_{w}(g)}{\|\cdot\|^2_{v}(g)}
    \right)
    =
    \sum_{k=0,1}(-1)^k
    \operatorname{Tr}\bigl(g(\beta+\beta^*) Q_{\overline{k},v}\bigr).
\end{equation}
\end{lemma}

Now we can state the variation formula of twisting element $H$.
\begin{thm}
    Let $X$ be a compact complex manifold, let $E$ be a holomorphic vector bundle over $X$, and let $H_v$ be a family of $\bpa$-closed $(0,\bline{1})$-forms such that the Dolbeault cohomology class $[H_v]$ is independent of $v$. Suppose that $G$ acts holomorphically on $X$ and $E$.
    
    By using $\epsilon_B$ to induce isomorphisms between $\det H^*(X,E,H_v)$, the twisted analytic torsions defined at different $v$ can be identified as metrics on $\det H^*(X,E,H_{v_0})$ at the base value $v_0$. Let $\beta:=B\wedge \cdot$, where $\bpa B=-\dfrac{\pa}{\pa v}H_v\bigg|_{v_0}$. Then we have:
    \begin{equation}
    \begin{aligned}
        &\frac{\pa}{\pa v}\log[\|\cdot\|_{E,H_v}(g)^2]\bigg|_{v_0}\\&=\sum_{k=0,1}(-1)^k \mathrm{Lim}_{t \to 0^+}\mathrm{Tr}(g(\beta+\beta^*) e^{-t\square_{\bline{k}}}),
        \end{aligned}
    \end{equation}
    where  $\mathrm{Lim}_{t \to 0^+}\mathrm{Tr}(g\beta e^{-t\square_{\bline{k}}})$ denotes the constant term of the asymptotic expansion of $\mathrm{Tr}(g\beta e^{-t\square_{\bline{k}}})$ at $t=0$.
\end{thm}
\subsection{Functorial properties}
\noindent In this section we present the basic functorial properties of equivariant twisted holomorphic analytic torsion. All of them can be obtained by the same method in \cite[Section 4]{MathaiWutwistedtorsion}. We first state several properties of twisted Dolbeault cohomology, which can be obtained by using Hodge decomposition.

\begin{prop}\label{prop:directsum}
    Let $X$ be a compact complex manifold, let $E_1,E_2$ be holomorphic vector bundles over $X$, and let $H$ be a $(0,\bline{1})$-form on $X$ satisfying $\bpa H=0$. Then, for every $k$, we have
    \begin{equation}
        H^{\bar{k}}(X,E_1\oplus E_2,H)\cong H^{\bar{k}}(X,E_1,H)\oplus H^{\bar{k}}(X,E_2,H).
    \end{equation}
\end{prop}

\begin{prop}[Serre duality]\label{prop:serre}
    Let $X$ be a compact complex manifold of complex dimension $n$, let $E$ be a holomorphic vector bundle over $X$, and let $H$ be a $(0,\bline{1})$-form on $X$ satisfying $\bpa H=0$. Let $K_X$ denote the canonical bundle of $X$. Then, for every $k$, we have
    \begin{equation}
        H^{\bar{k}}(X,E,H)\cong H^{\bline{n-k}}(X,E^*\otimes K_X,-H)^*.
    \end{equation}
\end{prop}

\begin{prop}[K\"{u}nneth formula]\label{prop:kunneth}
    Let $i=1,2$, and let $E_i$ be a holomorphic vector bundle over a compact complex manifold $X_i$. Let $H_i$ be a $(0,\bline{1})$-form on $X_i$ satisfying $\bpa H_i=0$, $i=1,2$. Let $\pi_i:X_1\times X_2 \to X_i$ be the natural projection. Let $E_1\boxtimes E_2=\pi_1^*E_1\otimes \pi_2^*E_2$ be the holomorphic vector bundle over $X_1\times X_2$, and let $H_1\boxplus H_2=\pi_1^*H_1+\pi_2^*H_2$ be the $\bpa$-closed $(0,\bline{1})$-form on $X_1\times X_2$. Then, for every $k$, we have
    \begin{equation}
        H^{\bar{k}}(X_1\times X_2,E_1\boxtimes E_2,H_1\boxplus H_2)\cong \bigoplus_{i+j=k} H^{\bar{i}}(X_1,E_1,H_1)\otimes H^{\bar{j}}(X_2,E_2,H_2).
    \end{equation}
\end{prop}

By the preceding three isomorphisms, we can now describe the functorial properties of twisted holomorphic analytic torsion.

\begin{prop}
    Under the hypotheses of Proposition~\ref{prop:directsum}. Suppose that $G$ acts holomorphically and isometrically on $X$, and that $E_1,E_2$ are equivariant holomorphic vector bundles over $X$. Then, for every $g\in G$, we have
    \begin{equation}
        \|\cdot\|_{E_1\oplus E_2,H}(g)^2= \|\cdot\|_{E_1,H}(g)^2\otimes \|\cdot\|_{E_2,H}(g)^2,
    \end{equation}
    where we use the canonical isomorphism
    \begin{equation*}
        \det H^*(X,E_1\oplus E_2,H)\cong \det H^*(X,E_1,H)\otimes \det H^*(X,E_2,H).
    \end{equation*}
\end{prop}

\begin{prop}
    Assume the conditions of Proposition~\ref{prop:serre}. Suppose that $G$ acts holomorphically and isometrically on $X$, and that $E$ is an equivariant holomorphic vector bundle over $X$. Then, for every $g\in G$, we have
    \begin{equation}
        \|\cdot\|_{E,H}(g)= \|\cdot\|_{E^*\otimes K_X,-H}(g^{-1})^{(-1)^{n+1}},
    \end{equation}
    where we use the canonical isomorphism
    \begin{equation*}    
     \det H^*(X,E,H)\cong \det H^*(X,E^*\otimes K_X,-H)^{(-1)^{n+1}}.
    \end{equation*}
\end{prop}

\begin{prop}
    Assume the conditions of Proposition~\ref{prop:kunneth}. Suppose that $G$ acts holomorphically and isometrically on $X_1$ and $X_2$, and that $E_1,E_2$ are equivariant holomorphic vector bundles over $X_1$ and $X_2$, respectively. Then, for every $g\in G$, we have
    \begin{equation}    
        \|\cdot\|_{E_1\boxtimes E_2,H_1\boxplus H_2}(g)= (\|\cdot\|_{E_1,H_1}(g))^{\otimes \chi_g(X_2,E_2,H_2)}\otimes (\|\cdot\|_{E_2,H_2}(g))^{\otimes \chi_g(X_1,E_1,H_1)},
    \end{equation}
    where we use the canonical isomorphism
    \begin{equation*}
    \begin{aligned}
        &\det H^\cdot(X_1\times X_2,E_1\boxtimes E_2,H_1\boxplus H_2)\cong \\ &(\det H^\cdot(X_1,E_1,H_1)^{\otimes \chi(X_2,E_2)})\otimes (\det H^\cdot(X_2,E_2,H_2))^{\otimes \chi(X_1,E_1)}.
    \end{aligned}
    \end{equation*}
    Here $\chi(X_i,E_i)$ denotes the holomorphic Euler characteristic of $E_i$ on $X_i$, and $\chi_g(X_i,E_i,H_i)$ denotes the equivariant holomorphic Euler characteristic of $(X_i,E_i,H_i)$ ,which is defined as $\chi_g(X_i,E_i,H_i)=\sum_{k=0}^{1}(-1)^k \mathrm{Tr}(g|_{H^{\bline{k}(X_i,E_i,H_i)}})$.
\end{prop}

\noindent{\bf AI use declaration.}
ChatGPT and other LLMs were used solely for language and grammar editing of this preprint. The author is a
non-native English speaker. The author
reviewed and edited the manuscript and takes
full responsibility for its content.

\bibliography{refs}

@book {MR1031992,
    AUTHOR = {Lawson, Jr., H. Blaine and Michelsohn, Marie-Louise},
     TITLE = {Spin geometry},
    SERIES = {Princeton Mathematical Series},
    VOLUME = {38},
 PUBLISHER = {Princeton University Press, Princeton, NJ},
      YEAR = {1989},
     PAGES = {xii+427},
      ISBN = {0-691-08542-0},
   MRCLASS = {53-02 (53A50 53C20 57R75 58G10)},
  MRNUMBER = {1031992},
MRREVIEWER = {N.\ J.\ Hitchin},
}

@book {Pseudodifferential,
    AUTHOR = {Shubin, M. A.},
     TITLE = {Pseudodifferential operators and spectral theory},
   EDITION = {Second},
      NOTE = {Translated from the 1978 Russian original by Stig I.
              Andersson},
 PUBLISHER = {Springer-Verlag, Berlin},
      YEAR = {2001},
     PAGES = {xii+288},
      ISBN = {3-540-41195-X},
   MRCLASS = {47G30 (35Sxx 58J40)},
  MRNUMBER = {1852334},
       DOI = {10.1007/978-3-642-56579-3},
       URL = {https://doi.org/10.1007/978-3-642-56579-3},
}

@article {MathaiWutwistedtorsion,
    AUTHOR = {Mathai, Varghese and Wu, Siye},
     TITLE = {Analytic torsion for twisted de {R}ham complexes},
   JOURNAL = {J. Differential Geom.},
  FJOURNAL = {Journal of Differential Geometry},
    VOLUME = {88},
      YEAR = {2011},
    NUMBER = {2},
     PAGES = {297--332},
      ISSN = {0022-040X,1945-743X},
   MRCLASS = {58J52 (58J10)},
  MRNUMBER = {2838268},
MRREVIEWER = {Rung-Tzung\ Huang},
       URL = {http://projecteuclid.org/euclid.jdg/1320067649},
}

@book {MR2093043,
    AUTHOR = {Huybrechts, Daniel},
     TITLE = {Complex geometry},
    SERIES = {Universitext},
      NOTE = {An introduction},
 PUBLISHER = {Springer-Verlag, Berlin},
      YEAR = {2005},
     PAGES = {xii+309},
      ISBN = {3-540-21290-6},
   MRCLASS = {32Qxx (14-01 32-01 32G05 32J25 53C55 53C56)},
  MRNUMBER = {2093043},
MRREVIEWER = {Richard\ P.\ Thomas},
}

@incollection {MathaiWuZ2gradedtorsion,
    AUTHOR = {Mathai, Varghese and Wu, Siye},
     TITLE = {Analytic torsion of {$\Bbb Z_2$}-graded elliptic complexes},
 BOOKTITLE = {Noncommutative geometry and global analysis},
    SERIES = {Contemp. Math.},
    VOLUME = {546},
     PAGES = {199--212},
 PUBLISHER = {Amer. Math. Soc., Providence, RI},
      YEAR = {2011},
      ISBN = {978-0-8218-4944-6},
   MRCLASS = {58J52 (57R56 58J10)},
  MRNUMBER = {2815136},
MRREVIEWER = {Rung-Tzung\ Huang},
       DOI = {10.1090/conm/546/10790},
       URL = {https://doi.org/10.1090/conm/546/10790},
}

@article {MR1262703,
    AUTHOR = {Bismut, J.-M. and Zhang, W.},
     TITLE = {Milnor and {R}ay-{S}inger metrics on the equivariant
              determinant of a flat vector bundle},
   JOURNAL = {Geom. Funct. Anal.},
  FJOURNAL = {Geometric and Functional Analysis},
    VOLUME = {4},
      YEAR = {1994},
    NUMBER = {2},
     PAGES = {136--212},
      ISSN = {1016-443X,1420-8970},
   MRCLASS = {58G26 (57Q10)},
  MRNUMBER = {1262703},
MRREVIEWER = {Friedbert\ Pr\"ufer},
       DOI = {10.1007/BF01895837},
       URL = {https://doi.org/10.1007/BF01895837},
}

@article {equnonresidue,
    AUTHOR = {Dave, Shantanu},
     TITLE = {An equivariant noncommutative residue},
   JOURNAL = {J. Noncommut. Geom.},
  FJOURNAL = {Journal of Noncommutative Geometry},
    VOLUME = {7},
      YEAR = {2013},
    NUMBER = {3},
     PAGES = {709--735},
      ISSN = {1661-6952,1661-6960},
   MRCLASS = {58J42 (19D55 58J40)},
  MRNUMBER = {3108693},
MRREVIEWER = {Yong\ Wang},
       DOI = {10.4171/JNCG/132},
       URL = {https://doi.org/10.4171/JNCG/132},
}

@article {Bismut2,
    AUTHOR = {Bismut, Jean-Michel},
     TITLE = {Equivariant immersions and {Q}uillen metrics},
   JOURNAL = {J. Differential Geom.},
  FJOURNAL = {Journal of Differential Geometry},
    VOLUME = {41},
      YEAR = {1995},
    NUMBER = {1},
     PAGES = {53--157},
      ISSN = {0022-040X,1945-743X},
   MRCLASS = {58G26 (32L10)},
  MRNUMBER = {1316553},
MRREVIEWER = {K.\ Dabrowski},
       URL = {http://projecteuclid.org/euclid.jdg/1214456007},
}

@book {Teschl,
    AUTHOR = {Teschl, Gerald},
     TITLE = {Mathematical methods in quantum mechanics},
    SERIES = {Graduate Studies in Mathematics},
    VOLUME = {99},
      NOTE = {With applications to Schr\"odinger operators},
 PUBLISHER = {American Mathematical Society, Providence, RI},
      YEAR = {2009},
     PAGES = {xiv+305},
      ISBN = {978-0-8218-4660-5},
   MRCLASS = {81-01 (47N50 81Qxx)},
  MRNUMBER = {2499016},
MRREVIEWER = {Rupert\ L.\ Frank},
       DOI = {10.1090/gsm/099},
       URL = {https://doi.org/10.1090/gsm/099},
}
\bibliographystyle{plain}

\end{document}